\documentclass{amsart}

\theoremstyle{definition}

\theoremstyle{remark}

\numberwithin{equation}{section}

\input{tcilatex}

\begin{document}
\title{ON INTEGRAL OPERATORS WITH OPERATOR VALUED KERNELS}
\author{Rishad Shahmurov}
\address{Department of Mathematics, Yeditepe University, Kayishdagi Caddesi,
34755 Kayishdagi, Istanbul, Turkey}
\email{shahmurov@hotmail.com}
\subjclass[2000]{Primary 45P05, 47G10, 46E40 }
\date{}
\keywords{Integral operators, Banach--valued Bochner spaces,
operator--valued kernels, interpolation of Banach spaces}

\begin{abstract}
Here $L_{q}\rightarrow L_{p}$ boundedness of integral operator with
operator-valued kernels is studied and the main result is applied to
convolution operators. Using these results $B_{q,r}^{s}\rightarrow
B_{p,r}^{s}$ regularity for Fourier multiplier operator is established.
\end{abstract}

\maketitle

\vspace{-7mm}

\section*{1. Introduction}

It is well-known that, solutions of inhomogeneous differential and integral
equations are represented by integral operators. In order to investigate
stability properties of these problems it is important to have boundedness
of corresponding integral operators in the studied function spaces. For
instance, the boundedness of Fourier multiplier operators play crucial role
in the theory of linear PDE especially in the study of maximal regularity
for elliptic and parabolic PDE. For exposition of the integral operators
with scalar valued kernels see $\left[ 3\right] $ and for the application of
multiplier theorems see $\left[ 2\right] .$

Maria Girardi and Lutz Weis $\left[ 4\right] $ recently proved that integral
operator%
\begin{equation*}
(Kf)(\cdot )=\dint\limits_{S}k(\cdot ,s)f(s)d\nu (s)\eqno(1)
\end{equation*}%
defines a bounded linear operator%
\begin{equation*}
K:L_{p}(S,X)\rightarrow L_{p}(T,Y)
\end{equation*}%
provided some measurability conditions and the following assumptions 
\begin{equation*}
\sup_{s\in S}\dint\limits_{T}\left\Vert k(t,s)x\right\Vert _{Y}d\mu (t)\leq
C_{1}\left\Vert x\right\Vert _{X}\text{ for all }x\in X,
\end{equation*}%
\begin{equation*}
\sup_{t\in T}\dint\limits_{S}\left\Vert k^{\ast }(t,s)y^{\ast }\right\Vert
_{X^{\ast }}d\nu (s)\leq C_{2}\left\Vert y^{\ast }\right\Vert _{Y^{\ast }}%
\text{ for all }y^{\ast }\in Y^{\ast }
\end{equation*}%
are satisfied. Inspired from $\left[ 4\right] $ we will show that $(1)$
defines a bounded linear operator%
\begin{equation*}
K:L_{q}(S,X)\rightarrow L_{p}(T,Y)
\end{equation*}%
if the kernel $k:T\times S\rightarrow B(X,Y)$ satisfies the following
conditions%
\begin{equation*}
\sup_{s\in S}\left( \int\limits_{T}\left\Vert k(t,s)x\right\Vert
_{Y}^{\theta }dt\right) ^{\frac{1}{\theta }}\leq C_{1}\left\Vert
x\right\Vert _{X}\text{ for all }x\in X,
\end{equation*}%
\begin{equation*}
\sup_{t\in T}\left( \int\limits_{S}\left\Vert k^{\ast }(t,s)y^{\ast
}\right\Vert _{_{X^{\ast }}}^{\theta }ds\right) ^{\frac{1}{\theta }}\leq
C_{2}\left\Vert y^{\ast }\right\Vert _{Y^{\ast }}\text{ for all }y^{\ast
}\in Y^{\ast }
\end{equation*}%
where 
\begin{equation*}
\frac{1}{q}-\frac{1}{p}=1-\frac{1}{\theta }
\end{equation*}%
for $1\leq q<\frac{\theta }{\theta -1}\leq \infty $ and $\theta \in \left[
1,\infty \right. ).$

Here $X$ and $Y$ are Banach spaces over the field $C$ and $X^{\ast }$ is the
dual space of $X.$ The space $B(X,Y)$ of bounded linear operators from $X$
to $Y$ is endowed with the usual uniform operator topology.

Now let us state some important notations from $\left[ 4\right] .$ A
subspace $Y$ of $X^{\ast }$ $\tau $-norms $X$, where $\tau \geq 1$, provided%
\begin{equation*}
\left\Vert x\right\Vert _{X}\leq \tau \sup_{x^{\ast }\in B(Y)}\left\vert
x^{\ast }(x)\right\vert \text{ }\forall x\in X.
\end{equation*}%
It is clear that if $Y$ $\tau $-norms $X$ then the canonical mapping 
\begin{equation*}
u:X\rightarrow Y^{\ast }\text{ with }\left\langle y,ux\right\rangle
=\left\langle x,y\right\rangle
\end{equation*}%
is an isomorphic embedding with 
\begin{equation*}
\frac{1}{\tau }\left\Vert x\right\Vert _{X}\leq \left\Vert u(x)\right\Vert
_{Y^{\ast }}\leq \left\Vert x\right\Vert _{X}.
\end{equation*}%
Let $\left( T,\sum_{T},\mu \right) $ and $\left( S,\sum_{S},\nu \right) $ be
positive measure spaces and 
\begin{equation*}
\sum_{S}^{\text{finite}}=\left\{ A\in \sum_{S}:\nu (A)<\infty \right\} ,%
\text{ }\sum_{S}^{\text{full}}=\left\{ A\in \sum_{S}:\nu (S\backslash
A)=0\right\} .
\end{equation*}%
${\LARGE \varepsilon }(S,X)$ will denote the space of finitely-valued and
finitely supported measurable functions from S into X i.e.%
\begin{equation*}
{\LARGE \varepsilon }(S,X)=\left\{
\sum\limits_{i=1}^{n}x_{i}1_{A_{i}}:x_{i}\in X,\text{ }A_{i}\in \sum_{S}^{%
\text{finite}},\text{ }n\in N\right\} .
\end{equation*}%
Note that ${\LARGE \varepsilon }(S,X)$ is norm dense in $L_{p}(S,X)$ for $%
1\leq p<\infty .$ Let $L_{\infty }^{0}(S,X)$ be the closure of ${\LARGE %
\varepsilon }(S,X)$ in $L_{\infty }(S,X)$ norm. In general $L_{\infty
}^{0}(S,X)\neq L_{\infty }(S,X)$ (see $\left[ 4,\text{ Proposition 2.2}%
\right] $ and $\left[ 4,\text{ Lemma 2.3}\right] )$ $.$

A vector valued function $f:S\rightarrow X$ is measurable if there is a
sequence $\left( f_{n}\right) _{n=1}^{\infty }\subset {\LARGE \varepsilon }%
(S,X)$ converging (in the sense of $X$ topology) to \ $f$ and it is $\sigma
(X,\Gamma )$-measurable provided $\left\langle f(\cdot ),x^{\ast
}\right\rangle :S\rightarrow K$ is measurable for each $x^{\ast }\in \Gamma
\subset X^{\ast }.$ Suppose $1\leq p\leq \infty $ and $\frac{1}{p}+\frac{1}{%
p^{\prime }}=1.$ There is a natural isometric embedding of $L_{p^{\prime
}}(T,Y^{\ast })$ into $\left[ L_{p}(T,Y)\right] ^{\ast }$ given by 
\begin{equation*}
\left\langle f,g\right\rangle =\int\limits_{T}\left\langle
f(t),g(t)\right\rangle d\mu (t)\text{ for }g\in L_{p^{\prime }}(T,Y^{\ast })%
\text{ and }f\in L_{p}(T,Y).
\end{equation*}%
Often $[E(X)]^{\ast }=E^{\ast }(X^{\ast })$, for example, provided $X$ has
the Radon-Nikodym property. Let us remind an important fact that if $X$ is
reflexive or if $X$ is separable, then $X$ has the Radon-Nikodym property.

\section*{2. $L_{q}\rightarrow L_{p}$ estimates for Integral Operators}

In this section we identify conditions on operator-valued kernel $k:T\times
S\rightarrow B(X,Y)$, extending theorems in [4] so that 
\begin{equation*}
\left\Vert K\right\Vert _{L_{q}(S,X)\rightarrow L_{p}(T,Y)}\leq C
\end{equation*}%
for $1\leq q\leq p$. To prove our main result we shall use interpolation
theorems of $L_{p}$ spaces. Therefore, we will study $L_{1}(S,X)\rightarrow
L_{\theta }(T,Y)$ and $L_{\theta ^{\prime }}(S,X)\rightarrow L_{\infty
}(T,Y) $ boundedness of integral operator (1). The following two conditions
are natural measurability assumptions on $k:T\times S\rightarrow B(X,Y)$.

\vspace{3mm}

\textbf{Condition 2.1.} For any $A\in \sum_{S}^{\text{finite}}$ and each $%
x\in X$

(a) there is $T_{A,x}\in \sum_{T}^{\text{full}}$ so that if $t\in T_{A,x}$
then the Bochner integral 
\begin{equation*}
\dint\limits_{A}k(t,s)xd\nu (s)\text{ exists},
\end{equation*}%
\ 

(b) $T_{A,x}:t\rightarrow \dint\limits_{A}k(t,s)xd\nu (s)$ defines a
measurable function from $T$ into $Y.$

Note that if $k$ satisfies the above condition then for each $f\in $ $%
{\LARGE \varepsilon }(S,X),$ there is $T_{f}\in \sum_{T}^{\text{full}}$ so
that the Bochner integral 
\begin{equation*}
\dint\limits_{S}k(t,s)f(s)d\nu (s)\text{ exists\ }
\end{equation*}%
and (1) defines a linear mapping%
\begin{equation*}
K:{\LARGE \varepsilon }(S,X)\rightarrow L_{0}(T,Y)
\end{equation*}%
where $L_{0}$ denotes the space of measurable functions.

\vspace{3mm}

\textbf{Condition 2.2.} The kernel $k:T\times S\rightarrow B(X,Y)$ satisfies
following properties:

(a) a real valued mapping $\left\Vert k(t,s)x\right\Vert _{X}^{\theta }$ is\
product measurable $\forall x\in X$

(b) there is $S_{x}\in \sum_{S}^{\text{full}}$ so that 
\begin{equation*}
\left\Vert k(t,s)x\right\Vert _{L_{\theta (T,Y)}}\leq C_{1}\left\Vert
x\right\Vert _{X}
\end{equation*}%
for $1\leq \theta <\infty $ and $x\in X.$

\vspace{3mm}

\textbf{Theorem 2.3.} Suppose $1\leq \theta <\infty $ and the kernel $%
k:T\times S\rightarrow B(X,Y)$ satisfies Condition 2.1 and Condition 2.2.
Then integral operator (1) acting on ${\LARGE \varepsilon }(S,X)$ extends to
a bounded linear operator%
\begin{equation*}
K:L_{1}(S,X)\rightarrow L_{\theta }(T,Y).
\end{equation*}

\textbf{Proof}. Let $f=\sum\limits_{i=1}^{n}x_{i}1_{A_{i}}(s)\in {\LARGE %
\varepsilon }(S,X)$ be fixed. Taking into account the fact that $1\leq
\theta ,$ using the general Minkowski-Jessen inequality and assumptions of
the theorem we obtain 
\begin{equation*}
\begin{array}{lll}
\Vert (Kf)(t)\Vert _{L_{\theta }(T,Y)} & \leq & \displaystyle\left[
\int\limits_{T}\left( \int\limits_{S}\Vert
k(t,s)\sum\limits_{i=1}^{n}x_{i}1_{A_{i}}(s)\Vert _{Y}d\nu (s)\right)
^{\theta }\,d\mu (t)\right] ^{\frac{1}{\theta }} \\ 
&  &  \\ 
\vspace{-3mm} & \leq & \displaystyle\int\limits_{S}\left(
\int\limits_{T}\Vert k(t,s)\sum\limits_{i=1}^{n}x_{i}1_{A_{i}}(s)\Vert
_{Y}^{\theta }d\mu (t)\right) ^{\frac{1}{\theta }}\,d\nu (s) \\ 
&  &  \\ 
& \leq & \displaystyle\int\limits_{S}\left[ \int\limits_{T}\left(
\sum\limits_{i=1}^{n}1_{A_{i}}(s)\Vert k(t,s)x_{i}\Vert _{Y}\right) ^{\theta
}d\mu (t)\right] ^{\frac{1}{\theta }}\,d\nu (s) \\ 
&  &  \\ 
& \leq & \displaystyle~\int\limits_{S}\sum\limits_{i=1}^{n}1_{A_{i}}(s)%
\left( \int\limits_{T}\Vert k(t,s)x_{i}\Vert _{Y}^{\theta }d\mu (t)\right) ^{%
\frac{1}{\theta }}\,d\nu (s)%
\end{array}%
\end{equation*}%
\begin{equation*}
\begin{array}{lll}
& \leq & \displaystyle\int\limits_{S}\sum\limits_{i=1}^{n}1_{A_{i}}(s)\Vert
k(t,s)x_{i}\Vert _{L_{\theta (T,Y)}}\,d\nu (s)\leq
C_{1}\sum\limits_{i=1}^{n}\Vert x_{i}\Vert
_{X}\int\limits_{S}1_{A_{i}}(s)\,d\nu (s) \\ 
& = & \displaystyle C_{1}\sum\limits_{i=1}^{n}\Vert x_{i}\Vert _{X}\nu
(A_{i})=C_{1}\left\Vert f\right\Vert _{L_{1}(S,X)}.%
\end{array}%
\end{equation*}%
Hence, $\left\Vert K\right\Vert _{L_{1}\rightarrow L_{\theta }}\leq C_{1}.$ %
\hbox{\vrule height7pt width5pt}

\vspace{3mm}

\textbf{Condition 2.4.} For each $y^{\ast }\in Z$ there is $T_{y^{\ast }}\in
\sum_{T}^{\text{full}}$ so that $\forall t\in T_{y^{\ast }}$

(a) a real valued mapping $\left\Vert k^{\ast }(t,s)x^{\ast }\right\Vert
_{X^{\ast }}^{\theta }$ is\ measurable $\forall x^{\ast }\in X^{\ast }$

(b) There is $S_{x}\in \sum_{S}^{\text{full}}$ so that 
\begin{equation*}
\left\Vert k^{\ast }(t,s)y^{\ast }\right\Vert _{L_{\theta (S,X^{\ast
})}}\leq C_{2}\left\Vert y^{\ast }\right\Vert _{Y^{\ast }}
\end{equation*}%
for $1\leq \theta <\infty $ and $x\in X.$ \vspace{3mm}

\textbf{Theorem 2.5.} Let $Z$ be a subspace of $Y^{\ast }$ that $\tau $%
-norms $Y.$ Suppose $1\leq \theta <\infty $ and $k:T\times S\rightarrow
B(X,Y)$ satisfies Condition 2.1 and Condition 2.4. Then integral operator
(1) acting on ${\LARGE \varepsilon }(S,X)$ extends to a bounded linear
operator%
\begin{equation*}
K:L_{\theta ^{\prime }}(S,X)\rightarrow L_{\infty }(T,Y).
\end{equation*}

\textbf{Proof}. Suppose $f\in {\LARGE \varepsilon }(S,X)$ and $y^{\ast }\in
Z $ are fixed. Let $T_{f\text{ }},$ $T_{y^{\ast }}\in \sum_{T}^{\text{full}}$
be corresponding sets due to Condition 2.1 and Condition 2.4. If $t\in T_{f%
\text{ }}\cap T_{y^{\ast }}$ then by using H\"{o}lder's inequality and
assumptions of the theorem we get 
\begin{equation*}
\begin{array}{lll}
|<y^{\ast },(Kf)(t)>_{Y}| & = & \displaystyle\left\vert \left\langle y^{\ast
},\int\limits_{S}k(t,s)f(s)d\nu \left( s\right) \right\rangle \right\vert \\ 
\vspace{-3mm} &  &  \\ 
& \leq & \int\limits_{S}\left\vert \left[ k^{\ast }(t,s)y^{\ast }\right]
f(s)\right\vert d\nu \left( s\right) \\ 
&  &  \\ 
& \leq & \displaystyle\Vert k^{\ast }(t,s)y^{\ast }\Vert _{L_{\theta
}(S,X^{\ast })}\Vert f(s)\Vert _{L_{\theta ^{\prime }}(S,X)} \\ 
&  &  \\ 
& \leq & \displaystyle C_{2}\Vert y^{\ast }\Vert \Vert f\Vert _{L_{\theta
^{\prime }}(S,X)}.%
\end{array}%
\end{equation*}

Since, $T_{f\text{ }}\cap T_{y^{\ast }}\in \sum_{T}^{\text{full}}$ and $Z$ $%
\tau $-norms $Y$%
\begin{equation*}
\left\Vert Kf\right\Vert _{L_{\infty }(T,Y)}\leq C_{2}\tau \Vert f\Vert
_{L_{\theta ^{\prime }}(S,X)}.
\end{equation*}%
Hence, $\left\Vert K\right\Vert _{L_{\theta ^{\prime }}\rightarrow L_{\infty
}}\leq \tau C_{2}.$ \hbox{\vrule height7pt width5pt}

In $\left[ 4,\text{ Lemma 3.9}\right] $ authors slightly improved
interpolation theorem $\left[ 1,\text{ Thm 5.1.2}\right] .$ The next lemma
is a more general form of $\left[ 4,\text{ Lemma 3.9}\right] .$

\vspace{3mm}

\textbf{Lemma 2.6.} Suppose a linear operator 
\begin{equation*}
K:{\LARGE \varepsilon }(S,X)\rightarrow L_{\theta }(T,Y)+L_{\infty }(T,Y)
\end{equation*}%
satisfies%
\begin{equation*}
\Vert Kf\Vert _{L_{\theta }(T,Y)}\leq C_{1}\left\Vert f\right\Vert
_{L_{1}(S,X)}\text{ and }\Vert Kf\Vert _{L_{\infty }(T,Y)}\leq
C_{2}\left\Vert f\right\Vert _{L_{\theta ^{\prime }}(S,X)}\eqno(2)
\end{equation*}%
Then, for $\frac{1}{q}-\frac{1}{p}=1-\frac{1}{\theta }$ and $1\leq q<\frac{%
\theta }{\theta -1}\leq \infty $ the mapping $K$ extends to a bounded linear
operator 
\begin{equation*}
K:L_{q}(S,X)\rightarrow L_{p}(T,Y)
\end{equation*}%
with%
\begin{equation*}
\left\Vert K\right\Vert _{L_{q}\rightarrow L_{p}}\leq \left( C_{1}\right) ^{%
\frac{\theta }{p}}\left( C_{2}\right) ^{1-\frac{\theta }{p}}.
\end{equation*}%
\textbf{Proof}. Let us first consider conditional expectation operator 
\begin{equation*}
(K_{0}f)=E\left( (Kf)_{1_{B}}|_{\sum }\right) \text{ }
\end{equation*}%
where $\sum $ is a $\sigma $-algebra of subsets of $B\in \sum_{T}^{\text{%
finite}}.$ From (2) it follows that%
\begin{eqnarray*}
\Vert K_{0}f\Vert _{L_{\theta }(T,Y)} &\leq &C_{1}\left\Vert f\right\Vert
_{L_{1}(S,X)}\text{ }<\infty \text{ and } \\
\Vert K_{0}f\Vert _{L_{\infty }(T,Y)} &\leq &C_{2}\left\Vert f\right\Vert
_{L_{\theta ^{\prime }}(S,X)}<\infty .
\end{eqnarray*}%
Hence, by Riesz-Thorin theorem ($\left[ 1,\text{ Thm 5.1.2}\right] $) we have%
\begin{equation*}
\Vert K_{0}f\Vert _{L_{p}(T,Y)}\leq \left( C_{1}\right) ^{\frac{\theta }{p}%
}\left( C_{2}\right) ^{1-\frac{\theta }{p}}\left\Vert f\right\Vert
_{L_{q}(S,X)}.\eqno(3)
\end{equation*}%
Now, taking\ into account (3) and using the same reasoning as in the proof
of $\left[ 4,\text{ Lemma 3.9}\right] $ one can easily show assertion of
this lemma. \hbox{\vrule height7pt width5pt}

\vspace{3mm}

\textbf{Theorem 2.7. (}\textit{Operator-valued Schur's test}\textbf{)} Let $%
Z $ be a subspace of $Y^{\ast }$ that $\tau $-norms $Y$ and $\frac{1}{q}-%
\frac{1}{p}=1-\frac{1}{\theta }$ for $1\leq q<\frac{\theta }{\theta -1}\leq
\infty .$ Suppose $k:T\times S\rightarrow B(X,Y)$ satisfies Condition 2.1,
Condition 2.2 and Condition 2.4 with respect to $Z$. Then integral operator
(1) extends to a bounded linear operator%
\begin{equation*}
K:L_{q}(S,X)\rightarrow L_{p}(T,Y)
\end{equation*}%
with%
\begin{equation*}
\left\Vert K\right\Vert _{L_{q}\rightarrow L_{p}}\leq \left( C_{1}\right) ^{%
\frac{\theta }{p}}\left( \tau C_{2}\right) ^{1-\frac{\theta }{p}}.
\end{equation*}

\textbf{Proof}. Combining Theorem 2.3, Theorem 2.5 and Lemma 2.6 we obtain
assertion of the theorem. \hbox{\vrule height7pt width5pt}

\vspace{3mm}

\textbf{Remark 2.8. }Note that choosing $\theta =1$ we get original results
in $\left[ 4\right] .$

For $L_{\infty }$ estimates (it is more delicate and based on ideas from
geometry Banach spaces) and weak continuity and duality results see $\left[ 4%
\right] .$ The next corollary plays important role in the Fourier Multiplier
theorems.

\vspace{3mm}

\textbf{Corollary 2.9.} Let $Z$ be a subspace of $Y^{\ast }$ that $\tau $%
-norms $Y$ and $\frac{1}{q}-\frac{1}{p}=1-\frac{1}{\theta }$ for $1\leq q<%
\frac{\theta }{\theta -1}\leq \infty .$ Suppose $k:R^{n}\rightarrow B(X,Y)$
is strongly measurable on $X,$ $k^{\ast }:R^{n}\rightarrow B(Y^{\ast
},X^{\ast })$ is strongly measurable on $Z$ and 
\begin{equation*}
\left\Vert kx\right\Vert _{L_{\theta }(R^{n},Y)}\leq C_{1}\left\Vert
x\right\Vert _{X}\text{ for all }x\in X,
\end{equation*}%
\begin{equation*}
\left\Vert k^{\ast }y^{\ast }\right\Vert _{_{L_{\theta }(R^{n},X^{\ast
})}}\leq C_{2}\left\Vert y^{\ast }\right\Vert _{Y^{\ast }}\text{ for all }%
y^{\ast }\in Y^{\ast }
\end{equation*}%
Then the convolution operator defined by 
\begin{equation*}
(Kf)(t)=\int\limits_{R^{n}}k(t-s)f(s)ds~\mbox{for}~t\in R^{n}
\end{equation*}%
satisfies $\Vert K\Vert _{L_{q}\rightarrow L_{p}}\leq \left( C_{1}\right) ^{%
\frac{\theta }{p}}\left( C_{2}\right) ^{1-\frac{\theta }{p}}.$

It is easy to see that $k:R^{n}\rightarrow B(X,Y)$ satisfies Condition 2.1,
Condition 2.2 and Condition 2.4 with respect to $Z.$ Thus, assertion of
corollary follows from Theorem 2.7.

\section*{3. Fourier Multipliers of Besov spaces}

In this section we shall indicate importance of Corollary 2.9 in the theory
of Fourier multipliers (FM). Thus we give definition and some basic
properties of operator valued FM and Besov spaces.

Consider some subsets $\{J_{k}\}_{k=0}^{\infty }$ and $\{I_{k}\}_{k=0}^{%
\infty }$ of $R^{n}$, where 
\begin{equation*}
J_{0}=\left\{ t\in R^{n}:|t|\leq 1\right\} ,~J_{k}=\left\{ t\in
R^{n}:2^{k-1}\leq |t|\leq 2^{k}\right\} ~\mbox{for}~k\in N
\end{equation*}%
and 
\begin{equation*}
I_{0}=\left\{ t\in R^{n}:|t|\leq 2\right\} ,~I_{k}=\left\{ t\in
R^{n}:2^{k-1}\leq |t|\leq 2^{k+1}\right\} ~\mbox{for}~k\in N.
\end{equation*}%
Let us define the partition of unity $\{\varphi _{k}\}_{k\in N_{0}}$ of
functions from $S(R^{n},R).$ Suppose $\psi \in S(R,R)$ is a nonnegative
function with support in $[2^{-1},2],$ which satisfies 
\begin{equation*}
\sum\limits_{k=-\infty }^{\infty }\psi (2^{-k}s)=1~\mbox{for}~s\in
R\backslash \{0\}
\end{equation*}%
and 
\begin{equation*}
\varphi _{k}(t)=\psi (2^{-k}|t|),~\varphi
_{0}(t)=1-\sum\limits_{k=1}^{\infty }\varphi _{k}(t)~\mbox{for}~t\in R^{n}.
\end{equation*}%
Let $1\leq q\leq r\leq \infty $ and $s\in R.$ The Besov space is the set of
all functions $f\in S^{\prime }(R^{n},X)$ for which 
\begin{equation*}
\begin{array}{lll}
\Vert f\Vert _{B_{q,r}^{s}(R^{n},X)}: & = & \displaystyle\left\Vert
2^{ks}\left\{ (\check{\varphi}_{k}\ast f)\right\} _{k=0}^{\infty
}\right\Vert _{l_{r}(L_{q}(R^{n},X))} \\ 
&  &  \\ 
& \equiv & \displaystyle\left\{ 
\begin{array}{ll}
\displaystyle\left[ \sum\limits_{k=0}^{\infty }2^{ksr}\Vert \check{\varphi}%
_{k}\ast f\Vert _{L_{q}(R^{n},X)}^{r}\right] ^{\frac{1}{r}} & \mbox{if}%
~r\neq \infty \\ 
&  \\ 
\displaystyle\sup\limits_{k\in N_{0}}\left[ 2^{ks}\Vert \check{\varphi}%
_{k}\ast f\Vert _{L_{q}(R^{n},X)}\right] & \mbox{if}~r=\infty%
\end{array}%
\right.%
\end{array}%
\end{equation*}%
is finite; here $q$ and $s$ are main and smoothness indexes respectively.
The Besov space has significant interpolation and embedding properties: 
\begin{equation*}
B_{q,r}^{s}(R^{n};X)=\left( L_{q}(R^{n};X),W_{q}^{m}(R^{d};X)\right) _{\frac{%
s}{m},r},
\end{equation*}%
\begin{equation*}
W_{q}^{l+1}(X)\hookrightarrow B_{q,r}^{s}(X)\hookrightarrow
W_{q}^{l}(X)\hookrightarrow L_{q}(X)\text{ where }l<s<l+1,
\end{equation*}%
\begin{equation*}
B_{\infty ,1}^{s}(X)\hookrightarrow C^{s}(X)\hookrightarrow B_{\infty
,\infty }^{s}(X)\text{ for }s\in \mathbf{Z},
\end{equation*}%
and%
\begin{equation*}
B_{p,1}^{\frac{d}{p}}(R^{d},X)\hookrightarrow L_{\infty }(R^{d},X)\text{ for 
}s\in \mathbf{Z},
\end{equation*}%
where $m\in N$ and $C^{s}(X)$ denotes the Holder-Zygmund spaces.

\vspace{3mm}

\textbf{Definition 3.1.} Let $X$ be a Banach space and $1\leq u\leq 2.$ We
say $X$ has Fourier type $u$ if 
\begin{equation*}
\Vert \mathcal{F}f\Vert _{L_{u^{\prime }}(R^{n},X)}\leq C\Vert f\Vert
_{L_{u}(R^{n},X)}~\mbox{for each}~f\in S(R^{N},X),
\end{equation*}%
where $\frac{1}{u}+\frac{1}{u^{\prime }}=1,~\mathcal{F}_{u,n}(X)$ is the
smallest $C\in \lbrack 0,\infty ].$ Let us list some important facts:

(i) Any Banach space has a Fourier type $1,$

(ii) $B$-convex Banach spaces have a nontrivial Fourier type,

(iii) Spaces having Fourier type $2$ should be isomorphic to a Hilbert
spaces.

From $\left[ 5,\text{ \textbf{Theorem 3.1}}\right] $ it follows the
following corollary:

\vspace{3mm}

\textbf{Corollary 3.2.} Let $X$ be a Banach space having Fourier type $u\in
\lbrack 1,2]$ and $1\leq \theta \leq u^{\prime }.$ Then the inverse Fourier
transform defines a bounded operator 
\begin{equation*}
\mathcal{F}^{-1}:B_{u,1}^{n\left( \frac{1}{\theta }-\frac{1}{u^{\prime }}%
\right) }(R^{n},X)\rightarrow L_{\theta }(R^{n},X).\eqno(4)
\end{equation*}%
\vspace{3mm}

\textbf{Definition 3.3.} Let $\left( E_{1}(R^{n},X)\text{, }%
E_{2}(R^{n},Y)\right) $ be one of the following systems, where $1\leq q\leq
p\leq \infty $ 
\begin{equation*}
(L_{q}(X),L_{p}(Y))\text{ }\mbox{or}~(B_{q,r}^{s}(X),B_{p,r}^{s}(Y)).
\end{equation*}%
A bounded measurable function $m:R^{n}\rightarrow B(X,Y)$ is called a
Fourier multiplier from $E_{1}(X)$ to $E_{2}(Y)$ if there is a bounded
linear operator 
\begin{equation*}
T_{m}:E_{1}(X)\rightarrow E_{2}(Y)
\end{equation*}%
such that 
\begin{equation*}
T_{m}(f)=\mathcal{F}^{-1}[m(\cdot )(\mathcal{F}f)(\cdot )]~\mbox{for each}%
\text{ }f\in S(X),\eqno(5)
\end{equation*}%
\begin{equation*}
T_{m}~\mbox{is}~\sigma (E_{1}(X),E_{1}^{\ast }(X^{\ast }))~\mbox{to}~\sigma
(E_{2}(Y),E_{2}^{\ast }(Y^{\ast }))~\mbox{continuous.}\eqno(6)
\end{equation*}%
The uniquely determined operator $T_{m}$ is the FM operator induced by $m.$
Note that if $T_{m}\in B(E_{1}(X),E_{2}(Y))$ and $T_{m}^{\ast }$ maps $%
E_{2}^{\ast }(Y^{\ast })$ into $E_{1}^{\ast }(X^{\ast })$ then $T_{m}$
satisfies the weak continuity condition (6).

For definition of Besov spaces and its basic properties we refer to [5].

Since(5) can be written in the convolution form%
\begin{equation*}
T_{m}(f)(t)=\int\limits_{R^{n}}\check{m}(t-s)f(s)ds\eqno(7)
\end{equation*}%
Corollary 2.9 and Corollary 3.2 can be applied to obtain $%
L_{q}(R^{n},X)\rightarrow L_{p}(R^{n},Y)$ regularity for (5).

\vspace{3mm}

\textbf{Theorem 3.4.} Let $X$ and $Y$ be Banach spaces having Fourier type $%
u\in \lbrack 1,2]$ and $p$, $q\in \left[ 1,\infty \right] $ so that $0\leq 
\frac{1}{q}-\frac{1}{p}\leq \frac{1}{u}.$ Then there is a constant $C$
depending only on $\mathcal{F}_{u,n}(X)$ and $\mathcal{F}_{u,n}(Y)$ so that
if 
\begin{equation*}
m\in B_{u,1}^{n\left( \frac{1}{u}+\frac{1}{p}-\frac{1}{q}\right)
}(R^{n},B(X,Y))
\end{equation*}%
then $m$ is a FM from $L_{q}(R^{n},X)$ to $L_{p}(R^{n},Y)$ with 
\begin{equation*}
\Vert T_{m}\Vert _{L_{q}(R^{n},X)\rightarrow L_{p}(R^{n},Y)}\leq CM_{u}(m)~%
\mbox{for each}
\end{equation*}%
where 
\begin{equation*}
M_{u}(m)=\inf \left\{ \Vert m(a\cdot )\Vert _{B_{u,1}^{n\left( \frac{1}{u}+%
\frac{1}{p}-\frac{1}{q}\right) }(R^{n},B(X,Y))}:a>0\right\} .
\end{equation*}

\textbf{Proof.} Let $\frac{1}{q}-\frac{1}{p}=1-\frac{1}{\theta }$ and $1\leq
q<\frac{\theta }{\theta -1}\leq \infty .$ Assume that $m\in S\left( B\left(
X,Y\right) \right) .$ Then $\check{m}\in S\left( B\left( X,Y\right) \right)
. $ Since $\mathcal{F}^{-1}\left[ m\left( a\cdot \right) x\right] \left(
s\right) =a^{-n}\check{m}\left( \frac{s}{a}\right) x,$~choosing appropriate $%
a$ and using $\left( 4\right) $ we obtain 
\begin{equation*}
\begin{array}{lll}
\left\Vert \check{m}x\right\Vert _{L_{\theta }\left( Y\right) } & = & %
\displaystyle\left\Vert \left[ m\left( a\cdot \right) x\right] ^{\vee
}\right\Vert _{L_{\theta }\left( Y\right) } \\ 
&  &  \\ 
& \leq & \displaystyle C_{1}\left\Vert m\left( a\cdot \right) \right\Vert
_{B_{u,1}^{n\left( \frac{1}{\theta }-\frac{1}{u^{\prime }}\right)
}}\left\Vert x\right\Vert _{X} \\ 
&  &  \\ 
& \leq & \displaystyle2C_{1}M_{u}(m)\left\Vert x\right\Vert _{X}%
\end{array}%
\eqno(8)
\end{equation*}%
where $C_{1}$ depends only on $\mathcal{F}_{u,n}(Y).$~Since $m\in S\left(
B\left( X,Y\right) \right) $ we have $\left[ m^{\ast }\right] ^{\vee }=\left[
\check{m}\right] ^{\ast }\in S\left( B\left( Y^{\ast },X^{\ast }\right)
\right) $ and $M_{u}(m)=M_{u}(m^{\ast }).$ Thus, in a similar manner as
above, we get

\begin{equation*}
\begin{array}{lll}
\left\Vert \left[ \check{m}\left( \cdot \right) \right] ^{\ast }y^{\ast
}\right\Vert _{L_{\theta }\left( Y\right) } & \leq & \displaystyle%
2C_{2}M_{u}(m)\left\Vert y^{\ast }\right\Vert _{Y^{\ast }}%
\end{array}%
\eqno(9)
\end{equation*}%
for some constant $C_{2}$ depending on $F_{u,n}(X^{\ast }).$ Hence by (8-9)
and Corollary 2.9%
\begin{equation*}
\left( T_{m}f\right) \left( t\right) )=\dint\limits_{R^{n}}\check{m}\left(
t-s\right) f\left( s\right) ds
\end{equation*}%
satisfies%
\begin{equation*}
\left\Vert T_{m}f\right\Vert _{L_{p}\left( R^{n},Y\right) }\leq
CM_{u}(m)\left\Vert f\right\Vert _{L_{q}\left( R^{n},X\right) }
\end{equation*}%
for all $p$, $q\in \left[ 1,\infty \right] $ so that $0\leq \frac{1}{q}-%
\frac{1}{p}\leq \frac{1}{u}.$ Now, taking into account the fact that $%
S\left( B\left( X,Y\right) \right) $ is continuously embedded in $%
B_{u,1}^{n\left( \frac{1}{u}+\frac{1}{p}-\frac{1}{q}\right) }(B(X,Y))$ and
using the same reasoning as $[5,$Theorem $4.3]$ one can easily prove the
general case $m\in B_{u,1}^{n\left( \frac{1}{u}+\frac{1}{p}-\frac{1}{q}%
\right) }$ and weak continuity of $T_{m}.$ \hbox{\vrule height7pt
width5pt}

\vspace{3mm}

\textbf{Theorem 3.5.} Let $X$ and $Y$ be Banach spaces having Fourier type $%
u\in \lbrack 1,2]$ and $p$, $q\in \left[ 1,\infty \right] $ are so that $%
0\leq \frac{1}{q}-\frac{1}{p}\leq \frac{1}{u}.$ Then, there exist a constant 
$C$ depending only on $\mathcal{F}_{u,n}(X)$ and $\mathcal{F}_{u,n}(Y)$ so
that if $m:R^{n}\rightarrow B(X,Y)$ satisfy 
\begin{equation*}
\varphi _{k}\cdot m\in B_{u,1}^{n\left( \frac{1}{u}+\frac{1}{p}-\frac{1}{q}%
\right) }(R^{n},B(X,Y))~\mbox{and}~M_{u}(\varphi _{k}\cdot m)\leq A\eqno(10)
\end{equation*}%
then $m$ is a FM from $B_{q,r}^{s}(R^{n},X)$ to $B_{p,r}^{s}(R^{n},Y)$ and $%
\Vert T_{m}\Vert _{B_{q,r}^{s}\rightarrow B_{p,r}^{s}}\leq CA$ for each $%
s\in R$ and $r\in \lbrack 1,\infty ].$

Taking into consideration the Theorem 3.4 one can easily prove above theorem
in a similar manner as $[5,$ Theorem $4.3].$

The following corollary provides a practical sufficient condition to check
(10).

\vspace{3mm}

\textbf{Lemma 3.6.} Let $n\left( \frac{1}{u}+\frac{1}{p}-\frac{1}{q}\right)
<l\in N$ and $\theta \in \lbrack u,\infty ]$. If $m\in C^{l}(R^{n},B(X,Y))$
and 
\begin{equation*}
\Vert D^{\alpha }m|_{I_{0}}\Vert _{L_{\theta }(B(X,Y))}\leq A,~\Vert
D^{\alpha }m_{k}|_{I_{1}}\Vert _{L_{\theta }(B(X,Y))}\leq A,\text{ }%
m_{k}(\cdot )=m(2^{k-1}\cdot ),
\end{equation*}%
for each $\alpha \in N^{n},~|\alpha |\leq l$ and $k\in N$ then $m$ satisfies
(10).

Using the fact that $W_{u}^{l}(R^{n},B(X,Y))\subset B_{u,1}^{n\left( \frac{1%
}{u}+\frac{1}{p}-\frac{1}{q}\right) }(R^{n},B(X,Y)),$ the above lemma can be
proven in similar fashion as $[5,$ Lemma 4.10$]$.

Choosing $\theta =\infty $ in the Lemma 3.6 we get the following corollary:

\vspace{3mm}

\textbf{Corollary 3.7. (}Mikhlin's condition) Let $X$ and $Y$ be Banach
spaces having Fourier type $u\in \lbrack 1,2]$ and $0\leq \frac{1}{q}-\frac{1%
}{p}\leq \frac{1}{u}.$ If $m\in C^{l}(R^{n},B(X,Y))$ satisfies 
\begin{equation*}
\left\Vert (1+|t|)^{|\alpha |}D^{\alpha }m(t)\right\Vert _{L_{\infty
}(R^{n},B(X,Y))}\leq A
\end{equation*}%
for each multi--index $\alpha $ with $|\alpha |\leq l=\left\lceil n\left( 
\frac{1}{u}+\frac{1}{p}-\frac{1}{q}\right) \right\rceil +1,$ then $m$ is a
FM from $B_{q,r}^{s}(R^{n},X)$ to $B_{p,r}^{s}(R^{n},Y)$ for each $s\in R$
and $r\in \lbrack 1,\infty ]$.

\vspace{3mm}

\textbf{Remark 3.8. }Corollary 2.13 particularly implies the following facts:

(a) If $X$ and $Y$ are arbitrary Banach spaces then $l=\left\lceil n\left( 
\frac{1}{p}+\frac{1}{q^{\prime }}\right) \right\rceil +1$,

(b) If $X$ and $Y$ be Banach spaces having Fourier type $u\in \lbrack 1,2]$
and $\frac{1}{q}-\frac{1}{p}=\frac{1}{u}$ then Miklin's condition with order
1 suffices for multiplier function to be a FM in $\left(
B_{q,r}^{s}(R^{n},X),B_{p,r}^{s}(R^{n},Y)\right) .$\textbf{\ }

(c\textbf{) }Let $X$ and $Y$ be Hilbert spaces and $p$, $q\in \left[
1,\infty \right] $ are so that $\frac{1}{q}-\frac{1}{p}=\frac{1}{2}.$ If $%
m\in C^{1}(R^{n},B(X,Y))$ satisfies%
\begin{equation*}
\sup_{t\in R^{n}}\left\Vert m(t)\right\Vert _{B(X,Y)}\leq A_{1},
\end{equation*}%
\begin{equation*}
\sup_{t\in R^{n}}\left\Vert (1+|t|)\frac{d}{dt}m(t)\right\Vert _{B(X,Y)}\leq
A_{2}
\end{equation*}%
then $m$ is a FM from $B_{q,r}^{s}(R^{n},X)$ to $B_{p,r}^{s}(R^{n},Y)$.


\begin{thebibliography}{9}
\bibitem{F} J. B\"{o}rgh, J. L\"{o}fstr\"{o}m, \textit{Interpolation spaces:
An introduction}, \textsl{Springer--Verlag}, Berlin 1976, Grundlehren der
Mathematischen Wissenschaften, No. 223.

\bibitem{H} R. Denk, M. Hieber, J. Pr\"{u}ss, $R$-boundedness, Fourier
multipliers and problems of elliptic and parabolic type, \textsl{Mem. Amer.
Math. Soc.} \textsl{166}(2003), p. 1-106.

\bibitem{B} G. Folland, Real analysis, Pure and Applied Mathematics, John
Wiley \& Sons Inc., New York, 1984, Modern techniques and their
applications, A \textsl{Wiley-Interscience Publication}.

\bibitem{C} M. Girardi, L. Weis, Integral operators with operator-valued
kernels,\textit{\ }\textsl{Journal of Mathematical Analysis and Applications}%
, \textsl{290 }(2004), no.1, pp.190-212.

\bibitem{E} M. Girardi, L. Weis, Operator--valued multiplier theorems on
Besov spaces, \textsl{Math. Nachr., \textbf{251 }}(2003), 34--51.
\end{thebibliography}
\end{document}